\documentclass[a4paper]{amsart}
\usepackage{amsmath,amssymb,amsfonts,graphicx,epsf}


\newtheorem{prop}{\ \ \ Proposition}

\theoremstyle{definition}

\theoremstyle{remark}

\let\texthat\^
\renewcommand{\r}[1]{(\ref{#1})}
\renewcommand\^{\widehat }

\newcommand{\Zz}{\mathcal{Z}}
\newcommand{\vphi}{\varphi}
\newcommand{\HH}{{\bf (H0)}}
\def\R{{\mathbf{R}}}
\newcommand{\frp}[2]{\frac{\partial #1}{\partial #2}}
\newcommand{\veps}{\varepsilon}
\def\F{\mathbf{F}}
\def\cn{\mathrm{cn\,}}
\def\sn{\mathrm{sn\,}}
\def\dn{\mathrm{dn\,}}
\def\Z{\mathbf{Z}}
\def\sign{\mathrm{sign}}

\begin{document}
 
\title[Sphere and cut locus at tangency points for 2-D ARS]{The sphere and the cut locus at a tangency point in two-dimensional almost-Riemannian geometry}
\author
{B.~Bonnard, G.~Charlot, R.~Ghezzi, G.~Janin}
\thanks
{G. Janin was supported by DGA/D4S/MRIS, under the supervision of J. Blanc-Talon, DGA/D4S/MRIS, 
Responsable de Domaine Ing\'enierie de l'Information. 
B. Bonnard and G. Charlot were supported by the ANR project GCM}



\begin{abstract}
We study the tangential case in 2-dimensional almost-Riemannian geometry.
We analyse the connection with the Martinet case in sub-Rieman\-nian geometry.
We compute estimations of the exponential map which allow us to describe the conjugate locus and the cut locus at a tangency point. We prove that this last one generically accumulates at the tangency point as an asymmetric cusp whose branches are separated by the singular set.
\end{abstract}

\maketitle


\section{Introduction}

\label{section1}
\def\D{\mathbf{D}}
In a series of recent papers [2,3,8],  2-dimensional almost-Riemannian geometry is investigated under generic conditions, giving rise to Gauss-Bonnet type results on compact oriented surfaces.

Roughly speaking, an almost-Riemannian structure (ARS for short) on an $n$-dimensional manifold $M$ can be defined locally by the data of $n$ vector fields playing the role of an orthonormal basis. Where the vector fields are linearly independent, they define a Riemannian metric. But the structure is richer along the set $\Zz$ where they are linearly dependent (see section \ref{section2} for a precise definition of ARS).

For 2-dimensional ARS, it was proven in [2] that generically the singular set $\Zz$ is an embedded submanifold of dimension 1 and only 3 types of points exist: the ordinary points where the metric is Riemannian, the Grushin points where the distribution $\Delta$ generated by the vector fields has dimension 1 and is transverse to $\Zz$, and the tangency points where $\Delta$ has dimension 1 and is tangent to $\Zz$.

The situation around ordinary and Grushin points is well known from the metric point of view, even if new considerations about curvature close to the Grushin points allow the authors to prove new results in [2,3,8]. 

These metrics have also been  studied in [5]. In that paper, the authors deduce a global model on the two-sphere of revolution
$S^{2}$ as a deformation of the round sphere,
the metric being
$$
g_{\lambda}=d\vphi^{2}+G_{\lambda}(\vphi)d\theta^{2},~\lambda\in[0,1],
$$
with $G_{\lambda}(X)=\frac{X}{1-\lambda X}$, where $X=\sin^{2}\vphi$, and
$(\vphi,\theta)$ are the spherical coordinates. In this representation, the singularity is
located at the equator: $\vphi=\pi/2$. This metric
appears in orbital transfer and, moreover,  the
homotopy is important to understand the behavior of the curvature.
In this
framework a short analysis tells us that for the generic model the symmetries
(of revolution and with respect to the equator) cannot be preserved and a non integrable model is obtained.

In this paper we analyse the situation at tangency points. The presence of these points is fundamental in the study of 2-dimensional ARS.

In [3], the authors provide a classification of oriented ARS on compact oriented surfaces in terms of the Euler number of the vector bundle corresponding to the structure (see \ref{section1bis} for definition) in presence of tangency points, generalizing a result of [2].
The construction of Gauss Bonnet type formulae is more intricated in presence of tangency points because of the geometry of the tubular neighborhoods of $\Zz$ close to the tangency points (see [3]).

It happens that the geometry close to tangency points is not well known and more intricated for many reasons. First, the computation of expansions of the wave front is more complicated and involves elliptic functions. Second, the nilpotent approximation is far from being generic as defined below. In particular the distribution of the nilpotent approximation is transversal to its singular set at the tangency point.

\begin{figure}[h]
\begin{center}
\includegraphics[scale=0.8]{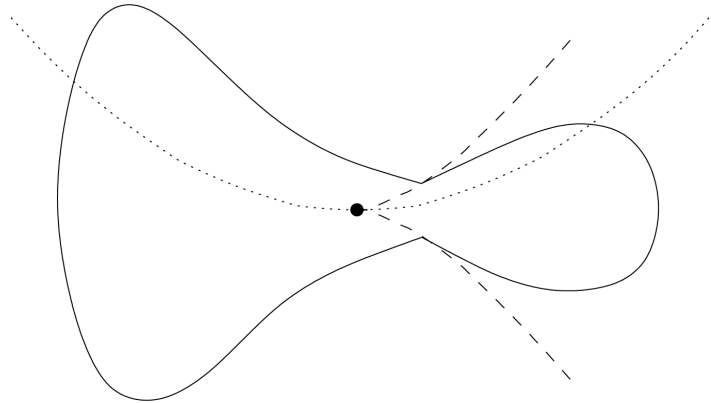}\label{cut}
\caption{{\small The sphere (solid line) and the cut locus (dashed line) at a tangency point in the generic case together with the singular set (dotted line)}}
\end{center}
\end{figure}

In this paper we focus on two points.
First, we analyse the connection between tangency points in 2-dimensional ARS and Martinet points in 3-dimensional sub-Riemannian structures. This allows us to obtain regularity properties of the distance function.
Second, we compute the jets of the exponential map which allows to estimate the conjugate locus and the cut locus at the tangency point. In particular we prove that, differing from the nilpotent case, the cut locus generically accumulates at the tangency point as an asymmetric cusp whose branches are locally separated by the singular set $\Zz$ (see figure \ref{cut}).

The paper is organised as follows. In section \ref{section1bis} we recall some basic definitions and results. In section \ref{section2} we show the relation between ARS and constant rank sub-Riemannian structures.  In section \ref{localanalysis} we analyse the special case of the nilpotent approximation as well as a generic model at a tangency point.
In section \ref{section3} we compute the asymptotic expansions of the exponential map at a tangency point. This allows us to estimate, in section \ref{section4}, the conjugate and cut loci at a tangency point, giving rise to a geometric interpretation of the first invariants in terms of the form of the cut locus.

\section{Basic definitions}
\label{section1bis}

An $n$-dimensional ARS is the data of a triple $(M,E,f)$ where $M$ is an n-dimensional manifold, E is a Euclidean bundle of rank $n$ over $M$ and $f$ is a morphism of vector bundles between $E$ and $TM$ preserving the basis $M$, such that the evaluation at any point $q \in M$ of the Lie algebra generated by
$\{f\circ \sigma \; |\; \sigma \mbox{ section of } E\}$ is $T_qM$.

From the control theory point of view, an ARS can be defined locally by the data of $n$ vector fields $(F_1,\dots, F_n)$ such that Lie$\{F_1,\dots,F_n\}_q=T_qM$ for all $q$. They define locally the following control dynamical system
\begin{equation}
\label{eq:PMP}
\dot q=\sum_{i=1}^n u_i F_i(q),\;\;\;\;\;\; \sum_{i=1}^n u_i^2= 1,
\end{equation}
the distance between two points $q_0$ and $q_1$ being by definition the minimal time needed to join $q_1$ from $q_0$ with this control system. We  also define the submodule $\Delta$ of the module of vector fields on $M$ generated locally by $(F_1,\dots, F_n)$ and the flag $\Delta_k$ by $\Delta_1=\Delta$, $\Delta_{k+1}=\Delta_k+[\Delta,\Delta_k]$.

In the following we deal with 2-dimensional ARS. Let us recall the following result proved in [2].

\begin{prop}
The following properties, denoted by \HH, are generic for 2-dimensional ARSs.
\begin{enumerate}
\item The singular set $\Zz$ is a one-dimensional embedded submanifold of $M$,
\item the points $q\in M$ where $\Delta_2(q)$ is one-dimensional are isolated,
\item $\Delta_3(q)=T_qM$ for all $q\in M$.
\end{enumerate}
Moreover, if a 2-dimensional ARS satisfies \HH, then for every point $q\in M$ there exist a neighborhood $U$ of $q$ and an orthonormal frame $(F_1,F_2)$ of the ARS on $U$ such that, up to a change of coordinates, $q=(0,0)$ and $(F_1,F_2)$ has one of the forms
$$
\begin{array}{lll}
(\F1) & F_1(x,y)=\frp{}{x}, & F_2(x,y)=e^{\phi(x,y)}\frp{}{y},\\
(\F2) & F_1(x,y)=\frp{}{x}, & F_2(x,y)=xe^{\phi(x,y)}\frp{}{y},\\
(\F3) & F_1(x,y)=\frp{}{x}, & F_2(x,y)=(y-x^2\psi(x))e^{\xi(x,y)}\frp{}{y},
\end{array}
$$
where $\phi$, $\psi$ and $\xi$ are smooth functions such that $\phi(0,y)=0$ and $\psi(0)>0$.
\end{prop}

\paragraph*{Remark}
In order to get the same notations as in [7], the normal form $(\F3)$ will be written in the following in  coordinates $(y,z)$
$$
\begin{array}{lll}
(\F3) & F_1(y,z)=\frp{}{y}, & F_2(y,z)=(z-y^2\psi(y))e^{\xi(y,z)}\frp{}{z}.
\end{array}
$$

For a 2-dimensional ARS satisfying \HH, we say that a point $q$ is
\begin{itemize}
\item ordinary if $\Delta(q)=T_qM$ (normal form $(\F1)$),
\item a Grushin point if the dimension of $\Delta(q)$ is one and $\Delta_2(q)=T_qM$ (normal form $(\F2)$),
\item a tangency point if the dimension of $\Delta(q)=\Delta_2(q)$ is one and $\Delta_3(q)=T_qM$ (normal form $(\F3)$).
\end{itemize}  

\medskip

In the normal form $(\F3)$, $(y,z)$ is a privileged coordinate system with weights respectively 1 and 3 (for definitions
of privileged coordinates and nilpotent approximation we refer the reader to [4]).

Consider the change of coordinates $\tilde y=y$ and $\tilde z=-z/(2\psi(0)e^{\xi(0)})$. Let us still denote $(\tilde y,\tilde z)$ by $(y,z)$. According to the weights  the jet up to order $0$ of the elements of the orthonormal frame in the normal form $(\F3)$ is
$$
\begin{array}{lll}
(\F3) & F_1(y,z)=\frp{}{y}, & F_2(y,z)=(\veps z+\frac{y^2}{2}+\veps'y^3+o_3(y,z))\frp{}{z}
\end{array}
$$
where $\veps= e^{\xi(0)}\neq0$, $\veps'=\frac{\psi'(0)+\psi(0)\frp{\xi}{y}(0)}{2\psi(0)}$ and $o_3(y,z)$ is a smooth function of order higher than 3 in the variables $y$ and $z$ with respect to their weights.

\section{Almost-Rieman\-nian geometry and sub-\-Rie\-man\-nian geometry}
\label{section2}
\subsection{Local desingularization of an $n$-dimensional ARS}\label{section2.1}
Let us present a classical construction.
Consider an $n$-dimensional ARS on $M$ and let $(F_1,\dots,F_n)$ be a local orthonormal frame on a neighborhood of $q$. Assume that $F_i(q)=0$ for $i>d$ where $d=\dim \Delta(q)$. Define
$$
\widetilde{M}=M\times \R^{n-d}=\{(x,y)|x\in M, y\in\R^{n-d}\}.
$$
Denote by $\pi_1$ and $\pi_2$ the canonical projections on the first and second factor of $\widetilde{M}$. Then we can define $\widetilde{F}_i$ by
$$
\begin{array}{ll}
{\pi_1}_*(\widetilde{F}_i)=F_i & \\
{\pi_2}_*(\widetilde{F}_i)=0 & \mbox{if } i\leq d, \\
{\pi_2}_*(\widetilde{F}_i)=\frac{\partial}{\partial y_{i-d}} & \mbox{if } i>d.
\end{array}
$$
Then the family $\{\widetilde{F}_1,\dots,\widetilde{F}_n\}$ has rank $n$ in the neighborhood of $q$ and defines the orthonormal frame of a sub-Riemannian metric on $\widetilde{M}$. Moreover, if $\Delta$ is bracket generating as a submodule of the Lie algebra of vector fields on $M$, then the same holds true for $\widetilde{\Delta}=\mbox{span}\{\widetilde{F}_1,\dots,\widetilde{F}_n\}$.
This metric on $\widetilde{M}$ is invariant with respect to translations in $\R^{n-d}$. Moreover, one can show that the curves between $q_0$ and $q_1$ minimizing the almost-Riemannian distance on $M$  are projections of the curves between $\{q_0\}\times \R^{n-d}$ and $\{q_1\}\times \R^{n-d}$ minimizing the sub-Riemannian distance on $\widetilde{M}$, a curve in $\widetilde{M}$ and its projection having the same length. This implies that the ball centered at $q$ of radius $r$ in $M$ is the projection of any sphere of radius $r$ centered at a point of the type $(q,y)$ in $\widetilde{M}$. Applying the Pontryagin Maximum Principle to the corresponding extremals in $\widetilde{M}$,  the transversality conditions to $\{q_0\}\times \R^{n-d}$ and to $\{q_1\}\times \R^{n-d}$ must be satisfied.

\subsection{Examples}

\subsubsection{The Grushin plane and the Heisenberg group}

The 2-dimensional ARS defined by the orthonormal frame
$\left\{F_1=\frp{}{x}, F_2=x\frp{}{y} \right\}$ in $\R^2$ is the Grushin plane.
It is the first example of almost-Riemannian structure with non empty singular locus, namely the $y$-axis.
Moreover, it is the nilpotent approximation at any Grushin point of a 2-dimensional ARS (for a precise definition of nilpotent approximation of a system we refer the reader to [4]).
If we apply the desingularization procedure, we find the sub-Riemannian metric defined by the orthonormal frame
$\left\{F_1=\frp{}{x},F_2=x\frp{}{y}+\frp{}{z} \right\}$ on $\R^3$, which is the Heisenberg metric.
The Heisenberg metric is the nilpotent approximation at any point of contact of a rank-2 sub-Riemannian structure defined on a 3-dimensional manifold, that is at any point $p$ where the rank-2 distribution satisfies $[\Delta,\Delta](p)=T_pM$. The Hamiltonian associated with the Grushin metric is
$$H_1=\frac{1}{2}(p_x^2+x^2p_y^2)$$
while the one related to the Heisenberg metric is
$$H_2=\frac{1}{2}(p_x^2+(xp_y+p_z)^2)$$
where $p_x$, $p_y$, $p_z$ are the dual coordinates to $x$, $y$ and $z$ in the cotangent bundle. Geodesics of the Heisenberg group projecting to geodesics of Grushin are those with $p_z=0$, respecting the transversality condition to the vertical lines given by the Pontryagin Maximum Principle.

In general, the relation between the cut and conjugate loci for the sub-\linebreak[1]Rie\-mannian metric on $\widetilde{M}$ and almost-Riemannian one on $M$ is not clear, the projection $\pi_1$ introducing singularities. The Grushin plane is a good illustration of this fact. The two geodesics for the Grushin metric starting at $(0,0)$ with the initial covectors $(p_x=1,p_y)$ and $(p_x=-1,p_y)$ are
\begin{eqnarray*}
x(t)&=&p_x \frac{\sin(p_y t)}{p_y},\\
y(t)&=&\frac{2p_y t-\sin(2p_y t)}{(2p_y)^2}.
\end{eqnarray*}
Hence, these two geodesics first intersect for $p_y\bar{t}=\pi$ and one can prove that $\bar t$ corresponds  to the cut time along them. Moreover, computing the Jacobian of the exponential mapping, one proves that the conjugate time $\tilde t$ satisfies $p_y \tilde t= \tan(p_y \tilde t)$. Lifting to the corresponding geodesics in the Heisenberg space starting at $(0,0,0)$ with initial covector $(p_x,p_y,p_z=0)$, one finds as third coordinate
\begin{eqnarray}
z(t)=p_x \frac{1-\cos(p_y t)}{p_y}.\nonumber
\end{eqnarray}
Hence, the two lifted geodesics do not intersect anymore at $\bar t$ and the computation of the Jacobian of the exponential mapping shows that the conjugate time for both lifted curves satisfies $\frac{p_y \tilde t}{2}= \tan(\frac{p_y \tilde t}{2})$. It corresponds to the second conjugate time in the Heisenberg case.

\subsubsection{The nilpotent approximation at a tangency point and the Martinet flat case}\label{TangencyMartinet}

Consider the normal form ($\F3$) at a tangency point as presented in section \ref{section1bis}. Recall that the weight of the variable $y$ is 1 and the weight of $z$ is 3. Hence the nilpotent approximation can be given in $(y,z)$ coordinates by the orthonormal frame $\left\{\frp{}{y}, \frac{y^2}{2}\frp{}{z} \right\}$, corresponding to the metric $g=dy^2+(\frac{y^2}{2})^{-2}dz^2$ and the Hamiltonian
$$H=\frac{1}{2}(p_y^2+\frac{y^4}{4}p_z^2).$$

Applying the desingularization procedure, one finds the orthonormal frame  
$$
\left\{\frp{}{y},~\frac{y^2}{2}\frp{}{z}+\frp{}{x} \right\}
$$
in $(y,z,x)$ coordinates in $\R^3$, whose corresponding Hamiltonian is
$$H=\frac{1}{2}(p_y^2+(p_x+\frac{y^2}{2}p_z)^2).$$
This lifted structure is the Martinet flat case of sub-Riemannian geometry. It is the nilpotent approximation at any Martinet point of a rank-2 sub-Riemannian structure defined on a 3-dimensional manifold, that is at any point $p$ where the rank-2 distribution satisfies $[\Delta,\Delta](p)=\Delta(p)$ and $[[\Delta,\Delta],\Delta](p)=T_pM$. The set of these points is called Martinet surface. Being the nilpotent approximation, the Martinet flat case will provide the starting point to analyse the general tangential case and it will allow to make some preliminary estimates of the sphere and the distance function using previous computations as in [7].

\section{Local analysis at a tangency point}\label{localanalysis}

In this section we focus on the following models in order to study the local situation  around tangency points for a generic 2-dimensional ARS.

\begin{enumerate}
\item {\em The nilpotent approximation (of order -1)}
$$
g_{-1}=dy^{2}+(\frac{y^{2}}{2})^{-2}dz^{2}
$$

\item {\em The generic model of order 0
$$
g_{0}=dy^{2}+(\veps z+\frac{y^{2}}{2}+\veps'y^{3})^{-2}dz^{2}
$$
where $\veps=\veps'=0$ gives the nilpotent approximation.}
\end{enumerate}

\subsection{ Analysis of the nilpotent model}\label{nil}
In this case, the desingularization procedure gives the orthonormal frame on $\R^3$
$$
F_{1}=\frp{}{x}+\frac{y^{2}}{2}\frp{}{z},~F_{2}=\frp{}{y}
$$
which generates the distribution
$$
\Delta=\ker(dz-\frac{y^{2}}{2}dx).
$$

\begin{prop}\label{prop-nilpotent}
Consider  the almost-Riemannian metric $g=dy^{2}+\frac{4}{y^{4}}dz^{2}$ on $\R^2$. The $y$-axis is the union of the two  geodesics starting at the origin with initial covectors $(p_y=\pm 1,p_z=0)$. The geodesics with initial covector $(p_y=\pm 1,p_z=\lambda\neq 0)$ are given by
\begin{eqnarray*}
y(t)&=&-p_y\frac{\sqrt{2}}{\sqrt{|\lambda|}}\cn(K+t\sqrt{|\lambda|}),\\
z(t)&=&\frac{\sign(\lambda)}{3|\lambda|^{3/2}}[t\sqrt{|\lambda|}+2\sn(K+t\sqrt{|\lambda|})\cn(K+t\sqrt{|\lambda|})\dn(K+t\sqrt{|\lambda|})],
\end{eqnarray*}
where $K$ is the complete elliptic integral of the first kind 
$$
\int_0^{\pi/2}\frac{d\varphi}{\sqrt{1-1/2\sin^2\varphi}}
$$
and $\cn,\sn,\dn$ denote the Jacobi elliptic functions of modulus $k=\frac{1}{\sqrt 2}$.
Moreover the following properties hold true.
\begin{enumerate}
\item The almost-Riemannian spheres centered at the origin are subanalytic. \label{subanalytic}
\item For $\lambda \neq 0$, the cut point coincides with  the first return to the $z$-axis that occurs at $t=2K/\sqrt{|\lambda|}$, where  two extremals
with the same length intersect. The cut locus from the origin is $\{(y,z)\mid y=0\}\setminus\{(0,0)\}$.
\item For $\lambda\neq 0$, the conjugate point corresponds to  $t\sim 3K/\sqrt{|\lambda|}$. The conjugate locus from the origin accumulates at the origin as a  set of the form
$\{(y,z)\mid z=\alpha y^3\}\cup\{(y,z)\mid z=-\alpha y^3\}\setminus\{(0,0)\}$, with  $\alpha \neq 0$.
\end{enumerate}
\end{prop}

{\bf Proof.}
Define $F_{3}=\frp{}{z}$ and $P_{i}=< p,F_{i}(q)>$, $i=1,2,3$. Using the Pontryagin Maximum Principle, the equations for normal extremals of the sub-Riemannian structure are given by the Hamiltonian system associated with $H=\frac{1}{2}(P_1^2+P_2^2)$, i.e.
$$
\begin{array}{rclrcl}
\dot{x} & = & P_{1}, & \dot{P}_{1} & = & yP_{2}P_{3},\\
\dot{y} & = & P_{2}, & \dot{P}_{2} & = & -yP_{1}P_{3},\\
\dot{z} & = & \frac{y^{2}}{2}P_{1}, & \dot{P}_{3} & = & 0.\end{array}
$$
There are three first integrals, namely $p_x,p_z=\lambda,H$. The normalization condition $H=1/2$ at $t=0$ gives
$$
P_{1}(0)^2+P_2(0)^2=1,
$$
hence, we
set $P_{1}(0)=\sin\vphi$, $P_{2}(0)=\cos\vphi$. The set of extremals
is invariant under the action of the group generated by the diffeomorphisms $(x,y,z)\mapsto(x,-y,z)$
and $(x,y,z)\mapsto(-x,y,-z)$. Therefore, it is sufficient to integrate the system with initial point $(0,0,0)$ and covector  $(\sin\vphi,\cos\vphi,\lambda)$ with $\lambda\geq 0$ and $\cos\vphi\geq 0$. Recall that the normal extremals for the almost-Riemannian structures are projections on the $(y,z)$ coordinates of the geodesics for the sub-Riemannian metric satisfying the transversality condition $p_x=0$. Hence, for  $\lambda=0$ we get $y(t)=\cos(\vphi) t, z(t)\equiv 0$.
Assume now $\lambda >0$ and set $k,k'$ such that $k^2=\frac{1-\sin\vphi}{2}$, $0<k,k'<1$ and  $k^2+{k'}^2=1$. Then we find
\begin{eqnarray*}
\dot{y}^{2}=(1-P_{1})(1+P_{1})&=&(1-p_{x}-\frac{y^{2}}{2}p_{z})(1+p_{x}+\frac{y^{2}}{2}p_{z})\\
                             &=&(2k^{2}-\frac{y^{2}}{2}\lambda)(2{k'}^{2}+\frac{y^{2}}{2}\lambda).
\end{eqnarray*}
Setting ${\mathrm y}(t)=\frac{\sqrt{\lambda}}{2k}y(t)$, the evolution equation for $\mathrm{y}$ is
$$
\frac{\dot{\mathrm{y}}^{2}}{\lambda}=(1-{\mathrm y}^{2})({k'}^{2}+k^{2}{\mathrm y}^{2}),
$$
that can be integrated, with $\dot{\mathrm{y}}(0)>0$, as $\mathrm{y}(t)=-\cn(K(k)+t\sqrt{\lambda},k)$, where
$$
K(k)=\int_0^{\pi/2}\frac{d\varphi}{\sqrt{1-k^2\sin^2\varphi}}.
$$
Hence
$$
y(t)=-\frac{2k}{\sqrt{\lambda}}~\cn(K(k)+t\sqrt{\lambda},k).
$$
Remark that the extremals that project on geodesics for the ARS  satisfy the transversality condition $p_x=0$ which implies $k^2=1/2$. Thus the $y$ coordinate of the geodesic with $p_x=0,p_y=1,\lambda>0$ is
$$
y(t)=-\frac{\sqrt{2}}{\sqrt{\lambda}}~\cn(K+t\sqrt{\lambda}),
$$  
where, to simplify notations, we denote $K(\sqrt{2}/2)$ by $K$ and we omit the dependence of the Jacobi function $\cn$ on the modulus.
To compute the  $z$ coordinate along the same geodesic, we use the primitive $\int\cn^{4}(K+ u)~du=\frac{2}{3}[\frac{1}{2}u+\sn (K+u)~\cn (K+u)~\dn(K+u)]$ (where $k=\sqrt{2}/2$). Then
$$
z(t)=\frac{1}{3\lambda^{3/2}}[t\sqrt{\lambda}+2~\sn(K+t\sqrt{\lambda})~\cn(K+t\sqrt{\lambda})~\dn(K+t\sqrt{\lambda})].
$$
Using the symmetries of the system we find the required expressions for the geodesics starting at the origin for the almost-Riemannian metric.
In particular, $y$ and $z$ are
quasi-homogeneous with respective weights $1$ and $3$. The cut instant of a geodesic coincides with the first
return to $y=0$ that occurs at  $t=2K/\sqrt{\lambda}$, thus the cut locus is the $z$-axis. The conjugate time satisfies $t\sim 3K/\sqrt{\lambda}$, whence the conjugate locus can be approximated  by the parametric curve
$$
y=-\frac{\sqrt 2}{\lambda^{1/2}},~z=\frac{K}{\lambda^{3/2}}
$$
(for the detailed proof see [1]).
\hfill$\blacksquare$

Remark that for the desingularized structure, i.e.,  the sub-Riemannian Martinet flat case, the sub-analyticity of the sphere is lost in
the abnormal direction for which $k\rightarrow 1$. This does not arise  for the almost-Riemannian structure, since geodesics satisfy $k^{2}=1/2$.

Property \ref{subanalytic} of proposition \ref{prop-nilpotent} can be generalized to the generic tangential case using [11], see also the computations in section \ref{section4}.

\subsection{Analysis of the generic model of order 0}\label{gen}

The objective of this section is to lift the generic model of order 0 into a constant rank sub-Riemannian model in order to analyse the role of the invariants in the optimal dynamics. A geometric interpretation will be given in section \ref{section4} in terms of the form of the cut locus.

Recall that from [7] the sub-Riemannian Martinet model of order zero is normalized to 
$$
(1+\alpha y)^{2}dx^{2}+(1+\beta x+\gamma y)^{2}dy^{2},
$$
where the distribution has the standard Martinet form
$$
2dz=y^{2}dx.
$$
In this normal form the parameters $\alpha,\beta,\gamma$ are related
to the geometric properties of the sphere with small radius and appear
in the pendulum interpretation of the extremals. More precisely, for
$\beta=0$ the extremal system is integrable while if $\beta$ is
non zero we have dissipation. In the integrable case the important
parameter is $\alpha$ and if it is non zero the abnormal direction
is strict. The role of the parameter $\gamma$ is unimportant and
it can be absorbed by reparameterization.

Consider the almost-Riemannian metric on $\R^2$ given by the orthonormal frame
\begin{equation}\label{AR0}
F_1=(\veps z + y^2/2+\veps' y^3)\frp{}{z},   F_2=\frp{}{y},
\end{equation}
where $\veps\neq 0$.
This metric can be seen as the generic model of order $0$ for an ARS in a neighborhood of a tangency point (use the normal form ($\F3$) and  weights of  coordinates).  Next proposition gives a possible lifting of the  model of order $0$ at a tangency point for an almost-Riemannian metric, showing the relation with the model of order $0$ of a Martinet type distribution for a sub-Riemannian metric.
\begin{prop}
The generic model of order $0$ for an ARS in a neighborhood of a tangency point lifts into  the sub-Riemannian
Martinet model of order zero $\frac{dx^{2}}{(\veps(1+x))^{2}}+\frac{dy^{2}}{(1+2\veps'y+o(y))^{2}}$
on the distribution $2dz-y^{2}dx=0$.
\end{prop}

{\bf Proof.}
Applying the desingularization procedure (see section \ref{TangencyMartinet}) to the almost-Riemannian metric defined by \r{AR0}, we get the sub-Riemannian metric in $\R^3$ defined by the orthonormal frame, still denoted by  $F_1,F_2$,
$$
F_{1}=\frp{}{x}+(\veps z+\frac{y^2}{2}+\veps' y^{3})\frp{}{z},~F_{2}=\frp{}{y}.
$$
One gets
$$
[F_{1},F_{2}]=-y(1+3\veps' y)\frp{}{z},
~~~[[F_1,F_2],F_2]=(1+6\veps' y)\frp{}{z},
$$
$$
[[F_{1},F_{2}],F_{1}]=\veps y(1+3\veps' y)\frp{}{z},
$$
hence the Martinet surface is the set  $\{(x,y,z)\mid y=0\}$. Moreover the singular control in the Martinet surface is defined by
$$
u_1 \det(F_1,F_2,[[F_1,F_2],F_1])+u_2 \det(F_1,F_2,[[F_1,F_2],F_2])=0
$$
which implies $u_2=0$. The corresponding  trajectories are solutions of
$$
\dot x=u_1,~ \dot y=0, \dot z=u_1\veps z.
$$
In order to build a coordinate system $(\tilde x,\tilde y,\tilde z)$ in which the  distribution has the normal form  $D=\mathrm{ker}~\omega$,
$\omega=d\tilde z-\tilde y^{2}/2d\tilde x$, we normalize the singular flow to lines parallel to the $\tilde x$-axis and lying  in the Martinet surface.
We consider the diffeomorphism
\begin{eqnarray*}
\tilde x&=&\frac{e^{-\veps x}}{-\veps}-1,\\
\tilde y&=&y\sqrt{1+2\veps'y}=y+\veps'y^2+o(y^2),\\
\tilde z&=&ze^{-\veps x}.\\
\end{eqnarray*}
The orthonormal frame in the new coordinate system becomes
$$
F_1=-\veps(1+\tilde x)\frp{}{\tilde x}-\veps(1+\tilde x)\frac{\tilde y^2}{2}\frp{}{\tilde z},~~F_2=(1+2\veps'\tilde{y}+o(\tilde y))\frp{}{\tilde y}.
$$
Hence the distribution is in the normal form $d\tilde z=\frac{\tilde y^2}{2}d\tilde x$
and the metric is given by
$$
g=\frac{d\tilde x^2}{\veps^2(1+\tilde x )^{2}}+\frac{d\tilde{y}^2}{(1+2\veps'\tilde y+o(\tilde y))^{2}}.
$$
\hfill$\blacksquare$

Introducing $F_3=\frp{}{z}$ and $P_i=<p,F_i>$, the extremal flow is given by
\begin{eqnarray*}
\dot X &=& -\veps(1+X)P_1,\\
\dot Y &=& (1+2\veps'Y+o(Y))P_2,\\
\dot Z &=& -\veps(1+X)\frac{Y^2}{2}P_1,\\
\dot P_1 &=& -\veps(1+X)Y(1+2\veps'Y+o(Y))P_2P_3,\\
\dot P_2 &=& \veps(1+X)Y(1+2\veps'Y+o(Y))P_1P_3,\\
\dot P_3 &=& 0.
\end{eqnarray*}
Setting $P_3=\lambda$ and using the time parameter $\tau$ such that $d\tau=(1+X)dt$, we can write
\begin{eqnarray*}
\frac{dX}{d\tau} &=& -\veps P_1,\\
\frac{dY}{d\tau} &=& \frac{(1+2\veps'Y+o(Y))}{1+X}P_2,\\
\frac{dZ}{d\tau} &=& -\veps\frac{Y^2}{2}P_1,\\
\frac{dP_1}{d\tau} &=& -\lambda\veps Y(1+2\veps'Y+o(Y))P_2,\\
\frac{dP_2}{d\tau} &=& \lambda\veps Y(1+2\veps'Y+o(Y))P_1.\\
\end{eqnarray*}
Define $\theta$ in $\R/2\pi\Z$ by $P_1=\cos(\theta)$ and $P_2=\sin(\theta)$. It satisfies
$$
\frac{d\theta}{d\tau}=\lambda\veps Y(1+2\veps'Y+o(Y))
$$
and then
$$
\frac{d^2\theta}{d\tau^2}=\lambda\veps\frac{1+4\veps'Y+o(Y)}{1+X}\sin(\theta),
$$
which can be approximated by
$$
\frac{d^2\theta}{d\tau^2}=\lambda\veps(1-X+4\veps'Y+o(Y))\sin(\theta).
$$
According to [7], this corresponds to a dissipative pendulum, the non nullity of the
parameter $\veps'$ inducing a coupling with the $y$-coordinate.
Note that more computations are necessary
to get the sub-Riemannian Martinet metric in the normal form of order 0, leading
to a true dissipative pendulum equation with no coupling with the
$x$ and $y$ variables, see [7].




\section{Asymptotics of the wave front}
\label{section3}

In this section we use the techniques and results from [7], developped in the
sub-Riemannian Martinet case, to compute asymptotics of the front from the tangency point for the generic model of order $0$ for ARS.
Remark that  the higher order terms in the expansion of the elements of the orthonormal frame play no role
in the estimation of the front and, consequently, in the estimation of the cut and conjugate loci from the tangency point (see section \ref{section4}), as one can check easily.

\begin{prop} \label{proposition31}
Consider the ARS on $\R^2$ defined by the orthonormal frame given in \r{AR0}. The extremals satisfying initial condition
$$(y,z,p_y,p_z)|_{t=0}=(0,0,\pm 1,\lambda)$$
with $|\lambda|\sim +\infty$ can be expanded as
\begin{eqnarray*}
y(t) & = & \eta Y^{0}(t/\eta)+\eta^2Y^{1}(t/\eta)+o(\eta^2),\\
z(t) & = & \eta^3 Z^{0}(t/\eta)+\eta^4 Z^{1}(t/\eta)+o(\eta^4),
\end{eqnarray*}
where $\eta=\frac{1}{\sqrt{|\lambda|}}$,
\begin{eqnarray}\label{termininilp}
&&\dot{Y}^{0}=P_{Y}^{0},\,\,\, \qquad\qquad\dot{P}_{Y}^{0}=-\frac{(P_{Z}^{0})^{2}(Y^{0})^{3}}{2},\\
&&\dot{Z}^{0}=\frac{P_{Z}^{0}(Y^{0})^{4}}{4}, \qquad \dot{P}_{Z}^{0}=0,\nonumber
\end{eqnarray}
with initial condition $(Y^0,Z^0,P_Y^0,P_Z^0)|_{t=0}=(0,0,\pm 1,\pm1)$ and
\begin{eqnarray}\label{terminiuno}
\dot{Y}^{1}&=&P_{Y}^{1},\\
\dot{Z}^{1}&=&\frac{1}{4}P_Z^1(Y^{0})^4+P_{Z}^{0}((Y^{0})^{3}Y^{1}+\veps Z^{0}(Y^{0})^{2}+\veps'(Y^0)^5),\nonumber\\
\dot{P}_{Y}^{1}&=&-P_{Z}^{0}P_{Z}^{1}(Y^{0})^{3}
-(P_{Z}^{0})^{2}(\frac{3}{2}(Y^{0})^{2}Y^{1}+\veps Z^{0}Y^{0}+\frac{5}{2}\veps' (Y^0)^4),\nonumber\\
\dot{P}_{Z}^{1}&=&-\frac{1}{2}(P_{Z}^{0})^{2}\veps(Y^{0})^{2},\nonumber
\end{eqnarray}
with initial condition $(Y^1,Z^1,P_Y^1,P_Z^1)|_{t=0}=(0,0,0,0)$.
\end{prop}

\paragraph*{Remark}
Computations in this case are similar to the ones of the Martinet
sub-Riemannian case. System \r{terminiuno} represents a variational equation whose
integration is related to the second-order equation
$$
\ddot{Y}^{1}+(\frac{3}{2}{P_{Z}^{0}}^{2}{Y^{0}}^{2})Y^{1}=K(Y^{0})
$$
where $Y^{0}$ is a periodic elliptic function.

{\bf Proof.}
Recall that  $y, z$ have weight $1$ and $3$, respectively. In order to have the standard Darboux form of order $1$,
we fix the weight $0$ for $p_{y}$ and $-2$ for $p_{z}$.

The Hamiltonian is
$$
H=\frac{1}{2}(p_{z}^{2}(\veps z+y^{2}/2+\veps' y^3)^2+p_{y}^{2})
$$
and the extremal flow is
$$
\begin{array}{rclrcl}
\dot{y} & = & p_{y},& \dot{p}_{y} & = & -p_{z}^{2}(\veps z+y^{2}/2+\veps' y^3)(y+3\veps' y^2),\\
\dot{z} & = & p_{z}(\veps z+y^{2}/2+\veps' y^3)^{2},& \dot{p}_{z} & = & -p_{z}^{2}(\veps z+y^{2}/2+\veps' y^3)\veps.
\end{array}
$$
According to the weights we set
$$
\begin{array}{rclrcl}
y & = & \eta Y, & p_{y} & = & P_{Y},\\
z & = & \eta^{3}Z, & p_{z} & = & \frac{P_{Z}}{\eta^{2}},
\end{array}
$$
where $\eta$ is a parameter. The evolution equations for $(Y,Z,P_Y,P_Z)$ are
\begin{eqnarray*}
\dot{Y} & = & \frac{P_{Y}}{\eta},\\
\dot{Z} & = & P_{Z}(\frac{Y^{4}}{4\eta}+\veps ZY^{2}+\veps'Y^5+
{\veps'}^2\eta Y^6+2\veps\veps'\eta ZY^3+\veps^{2}\eta Z^{2}),\\
\dot{P}_{Y} & = & -P_{Z}^{2}(\frac{Y^{3}}{2\eta}+\veps ZY+\frac{5}{2}\veps'Y^4
+3\eta\veps'Y^2(\veps Z+\veps'Y^3)),\\
\dot{P}_{Z} & = & -P_{Z}^{2}\veps(\frac{Y^{2}}{2}+\veps\eta Z+\veps'\eta Y^3).
\end{eqnarray*}
Considering the expansions with respect to $\eta$
\begin{eqnarray*}
&&Y =  Y^{0}+\eta Y^{1}+o(\eta), \quad P_{Y}  =  P_{Y}^{0}+\eta P_{Y}^{1}+o(\eta),\\
&&Z =  Z^{0}+\eta Z^{1}+o(\eta), \quad P_{Z} =  P_{Z}^{0}+\eta P_{Z}^{1}+o(\eta),
\end{eqnarray*}
 by identification we find that the leading terms satisfy
\begin{equation}\label{terminizero}
\begin{array}{rclrcl}
\dot{Y}^{0} & = & \frac{P_{Y}^{0}}{\eta}, & \dot{P}_{Y}^{0} & = & -\frac{({P_{Z}^{0})}^{2}{(Y^{0})}^{3}}{2\eta},\\
\dot{Z}^{0} & = & \frac{P_{Z}^{0}{(Y^{0})}^{4}}{4\eta}, & \dot{P}_{Z}^{0} & = & 0.\end{array}
\end{equation}
In particular $P_{Z}^{0}$ is constant. Setting $\lambda=p_z(0)$, for $\lambda\neq 0$ we can fix $\eta=1/\sqrt{|\lambda|}$ and then $P_{Z}^0$ is normalized to 1 or -1.

Introducing the time parameter $s=t\sqrt{|\lambda|}$ the equations \r{terminizero} for the first-order
approximation become
\begin{equation}\label{nilpotentsystem}
\begin{array}{rclrcl}
\frac{dY^{0}}{ds} & = & P_{Y}^{0}, & \frac{dP_{Y}^{0}}{ds} & = & -\frac{(P_{Z}^{0})^{2}(Y^{0})^{3}}{2},\\
\frac{dZ^{0}}{ds} & = & \frac{P_{Z}^{0}(Y^{0})^{4}}{4}, & P_{Z}^{0} & \equiv & \pm 1.\end{array}
\end{equation}
System \r{nilpotentsystem} coincides with the Hamiltonian system for  the nilpotent model that has been integrated in Proposition \ref{prop-nilpotent}, using elliptic functions with modulus $k$ such that $k^{2}=1/2$. The solution is given in Proposition \ref{prop-nilpotent}.
\begin{eqnarray*}
Y^0(s)&=&-P^0_Y(0)\sqrt{2}\cn(K+s),\\
Z^0(s)&=&\frac{P_{Z}^{0}}{3}(s+2\sn(K+s)\cn(K+s)\dn(K+s)),\\
P^0_Y(s)&=&P^0_Y(0)+(P^0_Y(0))^3(-1+\sqrt 2 \dn(K+s)\sn(K+s)),\\
P^0_Z(s)&\equiv& \pm1.
\end{eqnarray*}
Using $s=t/\eta$, the system for $(Y,Z,P_Y,P_Z)$ becomes
\begin{eqnarray*}
\frac{dY}{ds} & = & P_{Y},\\
\frac{dZ}{ds} & = & P_{Z}(\frac{Y^{4}}{4}+\eta(\veps ZY^{2}+\veps'Y^5)
+\eta^{2}({\veps'}^2Y^6+2\veps\veps'ZY^3+\veps^2Z^{2})),\\
\frac{dP_{Y}}{ds} & = & -P_{Z}^{2}(\frac{Y^{3}}{2}+\eta(\veps ZY+\frac{5}{2}\veps'Y^4)
+3\eta^2\veps'Y^2(\veps Z+\veps'Y^3)),\\
\frac{dP_{Z}}{ds} & = & -P_{Z}^{2}\veps(\eta\frac{Y^{2}}{2}+\eta^{2}(\veps Z+\veps'Y^3)).
\end{eqnarray*}
Hence,  identifying terms of order $0$, one gets
\begin{eqnarray*}
\frac{dY^1}{ds}&=&P_{Y}^{1},\\
\frac{dZ^1}{ds}&=&\frac{1}{4}P_Z^1(Y^{0})^4+P_{Z}^{0}((Y^{0})^{3}Y^{1}+\veps Z^{0}(Y^{0})^{2}+\veps'(Y^0)^5),\\
\frac{dP^1_Y}{ds}&=&-P_{Z}^{0}P_{Z}^{1}(Y^{0})^{3}
-(P_{Z}^{0})^{2}(\frac{3}{2}(Y^{0})^{2}Y^{1}+\veps Z^{0}Y^{0}+\frac{5}{2}\veps' (Y^0)^4),\\
\frac{dP^1_Z}{ds}&=&-\frac{1}{2}{P_{Z}^{0}}^{2}\veps{Y^{0}}^{2}.
\end{eqnarray*}
\hfill$\blacksquare$

\section{Geometric estimates of the conjugate and cut loci}
\label{section4}

\subsection{The conjugate locus}

The following result gives a description of the conjugate locus from a tangency point of a 2-dimensional ARS.

\begin{prop} Consider an ARS on $\R^2$ defined by the orthonormal frame
$$
F_1=(\veps z + y^2/2+\veps' y^3+o_3(y,z))\frp{}{z},   F_2=\frp{}{y}.
$$
Then there exists a constant $\alpha \neq 0$ such that the conjugate locus from $(0,0)$  accumulates at $(0,0)$ as the set  
$$
\{(y,z)\mid  z=\alpha y^{3}\}\cup \{(y,z)\mid   z=-\alpha y^{3}\}\setminus\{(0,0)\}.
$$
\end{prop}

{\bf Proof.}
Applying Proposition  \ref{proposition31}, the exponential map at $(0,0)$ is given by
$$
(\eta,s)\mapsto (\eta Y^0(s)+o(\eta),\eta^3 Z^0(s)+o(\eta^3)),
$$
where $s=t \sqrt{p_z(0)}$, $\eta$ parametrizes the initial covector as $(p_y(0)=\pm 1,p_z(0)=P_Z^0/\eta^2)$, and
$$
\begin{array}{l}
Y^0(s)=-P_Y^0(0)\sqrt{2}\cn(K+s),\\
Z^0(s)=P_Z^0\frac{1}{3}(s+2\sn(K+s)\cn(K+s)\dn(K+s)).
\end{array}
$$
The conjugate time is the first zero of the Jacobian of the exponential map. The Jacobian is equal, up to a multiplicative constant, to
$$
\eta^3(Y^0 \frac{\mathrm{d} Z^0}{\mathrm{d} s}-3Z^0\frac{\mathrm{d} Y^0}{\mathrm{d} s})+o(\eta^3).
$$
It was proven in [7] that the function
$$j(s)=Y^0(s) \frac{\mathrm{d} Z^0}{\mathrm{d} s}-3Z^0\frac{\mathrm{d} Y^0}{\mathrm{d} s}
$$
has its first positive zero at $s=s_0 \sim 3K$ and that $j'(s_0)\neq 0$. Hence, the conjugate time is of the form $s_0+o(1)$ where $o(1)$ is a continuous map going to zero when $\eta$ goes to zero.

In terms of the singularity theory,
this computation proves that the exponential map of a general two-dimensional ARS can be seen as a small deformation of the exponential map of the nilpotent case. In the nilpotent case, the exponential map has only stable singularities (folds) corresponding to the first conjugate locus. Hence, in the general case, the first conjugate locus also corresponds to folds and accumulates at $(0,0)$ as the set $\{(y,z) \;\;|\;\; (z-\alpha y^3)(z+\alpha y^3)=0 \}$, where $\alpha\neq 0$, see Proposition \ref{prop-nilpotent}.
\hfill$\blacksquare$

\subsection{The cut locus}

\begin{figure}[htbp]
   \begin{center}
      \includegraphics[scale=0.5]{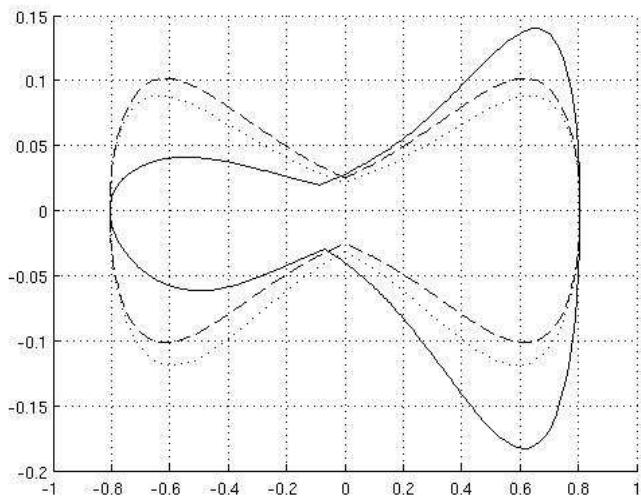}
   \end{center}
   \caption{The sphere of small radius for the nilpotent approximation (dotted line) and
for an example with $\epsilon'=0$ (dashed line) are symmetric. The two spheres
 are not $C^1$ at their intersection with the cut locus, which in both cases is the vertical axis. 
The sphere of small radius for the generic model of order zero in which $\varepsilon'\neq0$ 
(solid line) loses the symmetry. In this case, the cut locus is different from the
previous cases (see Proposition 6).}
   \label{spheres}
\end{figure}

In this section we provide a description of the cut locus at a tangency point for a generic ARS. As one can infer from the proof of the following proposition, the shape of the cut locus is determined only by the terms of order up to $0$ in the expansion of the elements of the orthonormal frame. Higher order terms do not contribute to the estimation of the way the cut locus approaches to the tangency point.
One can see in figure  \ref{spheres} small spheres for different values of $\veps$ and $\veps'$.

\begin{prop}\label{formacut} Consider the ARS on $\R^2$ defined by the orthonormal frame
$$
F_1=(\veps z + y^2/2+\veps' y^3+o_3(y,z))\frp{}{z},   F_2=\frp{}{y}.
$$
Then, if $\veps'\neq0$, there exist non zero constants $\alpha_1,\alpha_2$ such that the cut locus from $(0,0)$  accumulates at $(0,0)$ as the set  
$$
\{(y,z)\mid z>0,  z^{2}-\alpha_1 y^{3}=0\}\cup \{(y,z)\mid z<0,  z^{2}-\alpha_2 y^{3}=0\}.
$$
\end{prop}

{\bf Proof.}
In the following, we restrict our analysis to the upper half plane $z>0$, the computations being equivalent in the case $z<0$.
First of all, recall that for the nilpotent model ($\veps=\veps'=0$), the cut locus is $\{(y,z)\mid y=0\}\setminus \{(0,0)\}$ and the cut time on the geodesic with initial covector $(1,\lambda)$ is ${2K}/{\sqrt \lambda}$ which corresponds to the first intersection with the symmetric geodesic whose initial covector is  $(-1,\lambda)$, see Proposition \ref{prop-nilpotent}.  

Denote by $(Y^0,Z^0,P_Y^0,P_Z^0)$ the geodesic with initial condition $(0,0,1,1)$ for the ARS with $\veps=\veps'=0$, i.e.,
\begin{eqnarray*}
Y^{0}(s) & = & -\sqrt{2}\; \cn(K+s),\\
Z^{0}(s) & = & \frac{1}{3}(s+2\;\sn(K+s)\;\cn(K+s)\;\dn(K+s)),\\
P_Y^0(s) & = & \sqrt 2\; \dn(K+s)\;\sn(K+s),\\
P_Z^0 & \equiv & 1.
\end{eqnarray*}
Moreover, denote by  $Y^{1}, Z^{1}, P^{1}_Y, P^{1}_Z$ the terms of order $0$ in the expansion of the geodesic  for the ARS with
$\veps'=0$, i.e., solutions of the system \r{terminiuno}  with $\veps' =0$ and initial condition $(0,0,0,0)$. Set $({\mathcal Y}^1,{\mathcal Z}^1,{\mathcal P_Y}^1,{\mathcal P_Z}^1)$ to be
the terms of order $0$ in the expansion of the geodesic for the ARS with $\veps'\neq 0$, i.e., solutions of  system \r{terminiuno} with the initial condition $(0,0,0,0)$. Finally, define  four functions of $s$, $g_1, g_2, g_3$ and $g_4$, by
\begin{eqnarray*}
&&{\mathcal Y}^1=Y^1+\veps' g_1, \quad {\mathcal P_Y}^1={P_Y}^1+\veps' g_3,\\
&&{\mathcal Z}^1=Z^1+\veps' g_2, \quad {\mathcal P_Z}^1={P_Z}^1+\veps' g_4.
\end{eqnarray*}
Combining the equations satisfied by $Y^1$, $Z^1$, $P_Y^1$, $P_Z^1$ and by ${\mathcal Y}^1$, ${\mathcal Z}^1$, ${\mathcal P_Y}^1$, ${\mathcal P_Z}^1$, we  find that  $(g_1,g_2,g_3,g_4)$ satisfy the following system of ODEs
\begin{eqnarray}\label{equazioniperg}
&&\dot g_1 = g_3\qquad\qquad\qquad\qquad\,\dot g_3 = -\frac{3}{2}(Y^0)^2 g_1-\frac{5}{2}(Y^0)^4\nonumber\\
&&\dot g_2 = g_1 (Y^0)^3+(Y^0)^5\qquad g_4 \equiv 0,\nonumber
\end{eqnarray}
where $\dot g_i=dg_i/ds$,
with the initial conditions $
g_1(0)=g_2(0)=g_3(0)=0$.

Remark that if
$$(Y^0,Z^0,P_Y^0,P_Z^0,Y^1,Z^1,P_Y^1,P_Z^1,g_1,g_2,g_3)$$
is solution of (\r{termininilp},\r{terminiuno},\r{equazioniperg}) with the initial condition
$(0,0,1,1,0,0,0,0,0,0,0)$ then
$$(-Y^0,Z^0,-P_Y^0,P_Z^0,-Y^1,Z^1,-P_Y^1,P_Z^1,g_1,-g_2,g_3)$$
is also solution with the initial condition $(0,0,-1,1,0,0,0,0,0,0,0)$.
Moreover, one can compute numerically that  $g_1(2K)\sim-2\pi$, $g_2(2K)\sim-\pi$ and $g_3(2K)\sim 0$.

Let us compute the front at time $t=2K\eta_0$ close to the initial condition $\eta_0$,
that is for $\eta=\eta_0 + c \eta_0^2+o(\eta_0^2)$. Making Taylor expansions in terms of $\eta_0$,
one finds for the front corresponding to the initial conditions $p_y(0)=1$ and $p_z(0)={1}/{\eta^2}$
$$
\begin{array}{rcl}
y & = & \eta Y^0(\frac{2K\eta_0}{\eta})+\eta_0^2 (Y^1(2K)+\veps' g_1(2K))+o(\eta_0^2),\\
z & = & \eta^3 Z^0(\frac{2K\eta_0}{\eta}) + \eta_0^4(Z^1(2K)+\veps' g_2(2K))+o(\eta_0^4).
\end{array}
$$
Hence we obtain
$$
\begin{array}{rcl}
y & = & \eta_0^2 (Y^1(2K)+\veps' g_1(2K)+2Kc)+o(\eta_0^2),\\
z & = & \eta_0^3 Z^0(2K) + \eta_0^4(Z^1(2K)+\veps' g_2(2K)+ 3c Z^0(2K))+o(\eta_0^4),
\end{array}
$$
since $Y^0(2K)=\dot Z^0(2K)=0$ and $\dot Y^0(2K)=-1$.
For the front corresponding to the initial conditions $p_y(0)=-1$ and $p_z=1/\bar\eta^2$ where $\bar\eta=\eta_0 + c' \eta_0^2+o(\eta_0^2)$ one finds
$$
\begin{array}{rcl}
y & = & \eta_0^2 (-Y^1(2K)+\veps' g_1(2K)-2Kc')+o(\eta_0^2),\\
z & = & \eta_0^3 Z^0(2K) + \eta_0^4(Z^1(2K)-\veps' g_2(2K)+ 3c' Z^0(2K))+o(\eta_0^4).
\end{array}
$$
These expressions are affine with respect to parameters $c$ and $c'$, up to order $2$ for $y$ and $4$ for $z$ in the variable $\eta_0$. The two geodesics with initial covectors $p_y(0)=1,p_z(0)=1/\eta^2$ and $p_y(0)=-1,p_z(0)=1/\bar\eta^2$  intersect for
$$
\begin{array}{rcl}
c+c' & = & -\frac{Y^1(2K)}{K}+o(1),\\
c'-c & = & \frac{\veps' g_2(2K)}{K}+o(1)
\end{array}
$$
where $o(1)$ denotes any function going to 0 with $\eta_0$. Hence the intersection is for
$$
c=-\frac{\veps'g_2(2K)+Y^1(2K)}{2K}+o(1)
$$
and
$$
c'=\frac{\veps'g_2(2K)-Y^1(2K)}{2K}+o(1),
$$
which implies that the cut point is
$$
\begin{array}{rcl}
y_{\mathrm{cut}} & = & \eta_0^2\veps'(g_1(2K)-g_2(2K)) + o(\eta_0^2),\\
z_{\mathrm{cut}} & = & \eta_0^3 \frac{2K}{3}+o(\eta_0^3).
\end{array}
$$
Hence, if $\veps'\neq 0$, the upper branch of the cut locus from $(0,0)$ accumulates as the set $\{(y,z)\mid z>0, z^2=\alpha_1 y^3\}$, where
$$
\alpha_1 =\frac{4 K^2}{9 {\veps'}^3(g_1(2 K)-g_2(2K))^3}\sim -\frac{4 K^2}{9 {\veps'}^3\pi^3}.
$$
Similar computations show that the lower branch of the cut locus from $(0,0)$ accumulates as the set $\{(y,z)\mid z<0, z^2=\alpha_2 y^3\}$, where
$$
\alpha_2 =\frac{4 K^2}{9 {\veps'}^3(g_1(2 K)+g_2(2K))^3}\sim -\frac{4 K^2}{3^5 {\veps'}^3\pi^3}.
$$
\hfill$\blacksquare$

Remark that the case of generic ARS with $\veps'\neq 0$ is rather different from the sub-Riemannian Martinet case (see [7]), in which a similar argument cannot apply. Indeed, in the latter situation, the
asymptotic expansions for small time cannot be used, since there exists an abnormal
direction corresponding to the case where $k\rightarrow1$ and $K(k)\rightarrow\infty$.
Hence, we should use the asymptotic expansion for a time parameter tending to  $+\infty$,
which is clearly not valid.

\vskip+0.8cm

\centerline{\sc References} 

\vskip+0.4cm

1. A. Agrachev, B. Bonnard, M. Chyba, I. Kupka, Sub-Riemannian sphere in Martinet flat case,  {\em ESAIM/COCV},
\textbf{ 4}, 377--448, 1997.

2. A. Agrachev, U. Boscain and M. Sigalotti, A Gauss-Bonnet-Like
Formula on two-Dimensional Almost-Riemannian Manifolds, {\em Discrete
and Continuous Dynamical Systems}, vol. 20, No. \textbf{4}, April
2008, pp. 801-822.

3. A. Agrachev, U. Boscain G. Charlot, R. Ghezzi and M. Sigalotti,
Two-Dimensional
Almost-Riemannian Structures with Tangency Points. {\em Ann. Inst. H. Poincar\'e Anal. Non Lin\'eaire}. 27 (2010), pp. 793�-807.


4. A. Bella\"iche, The tangent space in sub-Riemannian geometry. {\em Dynamical systems, 3. J. Math. Sci.} (New York) 83 (1997),
no. 4, 461--476.

5. B. Bonnard, J.-B. Caillau, R. Sinclair and M. Tanaka,
Conjugate and cut loci of a two-sphere of revolution with application
to optimal control, {\em Ann. Inst. H. Poincar\'e Anal. Non Lin\'eaire}, 2009, vol. 26, no. 4, pp. 1081-1098.

6. B. Bonnard, J.-B. Caillau, Singular metrics on the
two-sphere in space mechanics, Preprint 2008, HAL, vol. 00319299, pp. 1-25

7. B. Bonnard and M. Chyba, M\'ethodes g\'eom\'etriques et analytiques
pour \'etudier l'application exponentielle, la sph\`ere et le front d'onde
en g\'eom\'etrie SR dans le cas Martinet, {\em ESAIM/COCV}, \textbf{4}, 245-334,
1999.

8. U. Boscain and M. Sigalotti, High-order angles in almost-Riemannian geometry, {\em Actes de S\'eminaires de Th\'eorie Spectrale et G\'eom\'etrie}, 25 (2008), pp 41-54.

9. J. Itoh, M. Tanaka, The Lipschitz continuity
of the distance function to the cut locus. {\em Trans. Amer. Math. Soc.}
353 (2001), no. 1, 21--40.

10. L.S. Pontryagin, V.G. Boltyanskii, R.V. Gamkrelidze, E.F. Mishchenko, The Mathematical
Theory of Optimal Processes, Interscience Publishers John Wiley and Sons, Inc, New York-London, 1962.

11. L. Rifford, E. Tr\'elat, Morse-Sard type results in sub-Riemannian geometry. {\em Math. Ann.} 332 (2005), no. 1,
145--159.

\ 

\ 

\centerline{(Received XXX)}

\ 

Authors' addresses:

B.~Bonnard, G.~Janin

Institut de Math{\'e}matiques de Bourgogne, 

UMR CNRS 5584,

BP 47870, 21078 Dijon Cedex, France

E-mail: bernard.bonnard{\string@}u-bourgogne.fr, gabriel.janin{\string@}u-bourgogne.fr
\vskip+2mm

G.~Charlot

Institut Fourier,

UMR CNRS 5582,

100 rue des Maths, BP 74, 38402 St Martin d'H{\`e}res, France

E-mail: Gregoire.Charlot{\string@}ujf-grenoble.fr
\vskip+2mm

R.~Ghezzi

SISSA/ISAS,

via Bonomea 265, 34136 Trieste, Italy

E-mail: ghezzi{\string@}sissa.it
\end{document}